\documentclass[12pt]{article}
\usepackage{amstex,cite}
\begin{document} 
\title{Permutation-type solutions to the Yang-Baxter and other
$n$-simplex equations}
\author{J. Hietarinta
\\ Department of Physics,
University of Turku\\ FIN-20014 Turku, Finland\\
E-mail: hietarin{@@}newton.tfy.utu.fi}
\date{\today}
\maketitle

\begin{abstract} We study permutation type solutions to
$n$-simplex equations, that is, solutions whose matrix form can be
written as $R_{i_i\dots i_n}^{j_1\dots j_n}=\prod_{\alpha=1}^{n}
\delta^{j_\alpha}_{A_{\alpha}^{\beta}i_{\beta} +B_\alpha}$ with some
$n\times n$ matrix $A$ and vector $B$, both over $\Bbb Z_D$.  With
this ansatz the $D^{n(n+1)}$ equations of the $n$-simplex equation
reduce to an $[\frac12n(n+1)+1]\times [\frac12n(n+1)+1]$ matrix
equation over $\Bbb Z_D$.  We have completely analyzed the 2-, 3- and
4-simplex equations in the generic $D$ case. The solutions show
interesting patterns that seem to continue to still higher simplex
equations.
\end{abstract}

\section{Introduction}
The Yang-Baxter equation (YBE, or 2-simplex equation) is the
fundamental equation of solvable models in (1+1)-dimensions. For
lattice models it guarantees the commutativity of the transfer matrix,
and for particle scattering it implies solvability through the
factorization of the scattering matrix \cite{YBL,tetL}. Therefore, in
order to construct interesting solvable models one needs interesting
solutions.  For this reason the YBE has been studied extensively and
indeed many solutions are known \cite{YBL,YBEs}, especially in the
two-state case \cite{complete}.

When one tries to generalize these solvable models to
(2+1)-dimensions, either by considering 3-dimensional lattices or the
scattering of straight strings, one obtains Zamolodchikov's
tetrahedron equation (3-simplex equation) as the fundamental equation
\cite{tetL,Label}, whose solutions are needed for further development.
Unfortunately only a few solutions are known for this equation
\cite{sim4,Kor,Label,some,Ser} and when one proceeds to still higher
dimensions and to the corresponding higher simplex equations very
little is known.
 
The difficulties associated with these equations come mainly from
sheer numbers, the $D$-state $n$-simplex equation is actually a set of
$D^{n(n+1)}$ equations on $D^{2n}$ variables (in the non-constant case
$(n+1)D^{2n}$ variables).  Because of this one is forced to make
rather restrictive ansatze in order to obtain any solutions at all.
One method is to take some definite high level structure (Lie algebra,
chiral Potts) coming from somewhere else and apply it to the present
situation.  Our approach is complementary to this, the ansatz given
below is defined in rather simple terms and we will then determine
{\em all} solutions within this class.

Let us recall the standard setup for the $n$-simplex equations.  As
usual we assume that we have linear operators ${\cal R}$ which for the
$n$-simplex case are assumed to act on a product of $n$ identical
vector spaces $V$, i.e., ${\cal R}:V^{\otimes n} \to V^{\otimes n}$.
Let $e_i$ be the $D$ basis vectors of $V$. Since we want to do algebra
with the indices of the basis vectors it would be nice if the indexing
formed a finite field. If $D$ is prime this is possible with $\Bbb
Z_D$, integers modulo $D$, which is what we consider in this paper.
However, some aspects of the following derivation works even if the
indices just form a ring, for example with $\Bbb Z_4$.

To the operator ${\cal R}$ we associate a numerical matrix with $n$
pairs of indices by
\begin{equation}
{\cal R}\,(e_{i_1}\otimes\dots\otimes e_{i_n})= R_{i_1\dots
i_n}^{j_1\dots j_n}\,(e_{j_1}\otimes\dots\otimes e_{j_n}).
\end{equation}
(Here and elsewhere in this paper summation over repeated indices is
assumed.)  The $n$-simplex equation itself is defined on $V^{\otimes
[n(n+1)/2]}$, and the linear operators operate trivially in all but
the $n$ spaces indicated by the subscripts, e.g.,$ {\cal R}_{12}\,
(e_{i_1}\otimes e_{i_2}\otimes e_{i_3})= R_{i_1i_2}^{j_1j_2}\,
(e_{j_1}\otimes e_{j_2}\otimes e_{i_3})$, or in the general case with
$K_\alpha\in\{1,\dots N\},\, N=\frac12n(n+1)$,
\begin{equation}
(R_{K_1\dots K_n})_{i_1\dots i_N}^{j_1\dots j_N}=
R_{i_{K_1}\dots i_{K_n}}^{j_{K_1}\dots j_{K_n}}{\hskip -0.3cm}
\prod\begin{Sb} k=1\\k\ne K_\alpha,\,\forall \alpha 
\end{Sb}^N{\hskip -0.3cm}\delta_{i_k}^{j_k}.
\label{E:Rmatd}
\end{equation}

In this paper we consider the first few constant simplex equations,
those given by the 2-simplex or vertex Yang-Baxter equation
\begin{equation}
{\cal R}_{12}{\cal R}_{13}{\cal R}_{23}=
{\cal R}_{23}{\cal R}_{13}{\cal R}_{12},
\label{E:2s}
\end{equation}
the 3-simplex or tetrahedron equation
\begin{equation}
{\cal R}_{123}{\cal R}_{145}{\cal R}_{246}{\cal R}_{356}=
{\cal R}_{356}{\cal R}_{246}{\cal R}_{145}{\cal R}_{123},
\label{E:3s}
\end{equation}
and the 4-simplex equation
\begin{equation}
{\cal R}_{1234}{\cal R}_{1567}{\cal R}_{2589}{\cal R}_{3680}{\cal
R}_{4790}= {\cal R}_{4790}{\cal R}_{3680}{\cal R}_{2589}{\cal
R}_{1567}{\cal R}_{1234}.
\label{E:4s}
\end{equation}
In terms of the multi-indexed matrices defined in (\ref{E:Rmatd})
the above operator equations imply, respectively,
\begin{equation}
R_{j_2j_3}^{k_2k_3} \,R_{j_1k_3}^{k_1l_3}\, R_{k_1k_2}^{l_1l_2}
=R_{j_1j_2}^{k_1k_2}\, R_{k_1j_3}^{l_1k_3}\, R_{k_2k_3}^{l_2l_3},
\end{equation}
\begin{equation}
R_{j_3j_5j_6}^{k_3k_5k_6}\,R_{j_2j_4k_6}^{k_2k_4l_6}
\,R_{j_1k_4k_5}^{k_1l_4l_5}\,R_{k_1k_2k_3}^{l_1l_2l_3}
=R_{j_1j_2j_3}^{k_1k_2k_3}\,R_{k_1j_4j_5}^{l_1k_4k_5}
\,R_{k_2k_4j_6}^{l_2l_4k_6\,}R_{k_3k_5k_6}^{l_3l_5l_6},
\end{equation}
\begin{align}
R_{j_4j_7j_9j_0}^{k_4k_7k_9k_0}
\,R_{j_3j_6j_8k_0}^{k_3k_6k_8l_0}\,&R_{j_2j_5k_8k_9}^{k_2k_5l_8l_9}
\,R_{j_1k_5k_6k_7}^{k_1l_5l_6l_7}\,R_{k_1k_2k_3k_4}^{l_1l_2l_3l_4}
=\nonumber\\
&R_{j_1j_2j_3j_4}^{k_1k_2k_3k_4}\,R_{k_1j_5j_6j_7}^{l_1k_5k_6k_7}
\,R_{k_2k_5j_8j_9}^{l_2l_5k_8k_9}\,R_{k_3k_6k_8j_0}^{l_3l_6l_8k_0}
\,R_{k_4k_7k_9k_0}^{l_4l_7l_9l_0}.
\end{align}
In addition to the above some other equations have appeared in the
literature, e.g., the Frenkel--Moore equation \cite{FM}. For a general
formulation of the various type of equations, see \cite{CS}.

\section{Formulation with the permutation ansatz} 
In this paper we consider only permutation type operators, that is
those which transform one product of basis vectors into another simple
product.  In the matrix form this means that there is precisely one
nonzero ($=1$) entry in each column and row.  For the $D$-state
$n$-simplex equation there are $(D^n)!$ different matrices to
consider, and a brute force check of them is out of question except for
$D=2,\, n=2$ which contains 24 permutation matrices (the next cases
$(2^3)!=40320$ and $(3^2)!=362880$ might still be possible).  We will
therefore make the further assumption that the dependence between the
basis vectors is {\em linear}, that is,
\begin{equation}
{\cal R}\,(e_{i_1}\otimes\dots\otimes e_{i_n})=
e_{A_1^\alpha i_\alpha+B_{1}}\otimes
\dots\otimes e_{A_n^\alpha i_\alpha+B_{n}},
\end{equation}
(where the summation over $\alpha$ runs from 1 to $n$) for some
nonsingular $n\times n$ matrix $A$ and $n$-vector $B$, both having
entries from $\Bbb Z_D$.  In terms of the $R$-matrix this means that
\begin{equation}
R_{i_1\dots i_n}^{j_1\dots j_n}=
\delta^{j_1}_{A_{1}^{\alpha}i_\alpha+B_1}\cdots
\delta^{j_n}_{A_{n}^{\alpha}i_\alpha+B_n}
\equiv
\delta(A,B).
\label{E:def1}
\end{equation}

The main advantage of this ansatz is that the problem of solving the
$D$-state $n$-simplex equation can be reduced to handling ordinary
matrices over $\Bbb Z_D$, as will be shown below.  This simplifies
the problem considerably.  Furthermore, although possible applications
normally imply further conditions on the solutions, permutation
matrices are such fundamental objects that there is a good change they
are acceptable in most cases, and we believe that the ansatz is not
an unnatural starting point.

In order to write the $n$-simplex equations in terms of $A$ and $B$
let us further define (in analogue with (\ref{E:Rmatd}))
\begin{eqnarray}
(A_{K_1\dots K_n})_i^j&=&\left\{
\begin{array}{ll}
A_{\alpha}^{\beta},&\mbox{ if }i=K_{\alpha},\,j=K_{\beta} 
\mbox{ for some }\alpha,\,\beta,\\
\delta_i^j, &\mbox{ otherwise,}
\end{array}\right.\label{E:Axxdef}\\
(B_{K_1\dots K_n})_i&=&\left\{
\begin{array}{ll}
B_i,&\mbox{ if }i=K_{\alpha},\,\mbox{ for some }\alpha,\\
0, &\mbox{ otherwise,}
\end{array}\right.
\end{eqnarray}
so that
\begin{equation}
(R_{K_1\dots K_n})_{i_1\dots i_N}^{j_1\dots j_N}=
\prod_{\mu=1}^N \delta^{j_n}_{(A_{K_1\dots K_n})_{\mu}^{\nu}i_\nu
+(B_{K_1\dots K_n})_\mu},
\end{equation}
where now the $\nu$ summation runs from 1 to $N$.

In the homogeneous case, that is with $B\equiv 0$, the above
correspondence between $R$ and $A$ means that the $n$-simplex equation
with ansatz (\ref{E:def1}) becomes an $N\times N$ matrix equation over
$\Bbb Z_D$.  For example, the 2-simplex equation becomes
\begin{equation}
(A_{12})_i^k(A_{13})_k^m(A_{23})_m^l=(A_{23})_i^k(A_{13})_k^m(A_{12})_m^l,
\label{E:Aeq}
\end{equation}
where $A_{KL}$ are $3\times3$ matrices with entries from $\Bbb Z_D$ as
given in (\ref{E:Axxdef}) (for the explicit form see (\ref{E:2sima})).

In the non-homogeneous case with $B\ne 0$ a matrix formulation can also be
obtained, if we add a fictitious index space 0 and write
\begin{equation}
R_{i_1\dots i_n}^{j_1\dots j_n}=
\delta^{j_1}_{A_{1}^{\alpha}i_\alpha+B_1i_0}\cdots
\delta^{j_n}_{A_{n}^{\alpha}i_\alpha+B_ni_0}\delta_{i_0}^{j_0}
=\prod_{\mu=0}^n\delta^{j_\mu}_{\tilde A_{\mu}^{\nu}i_\nu}.
\label{E:def10}
\end{equation}
When this is immersed in the larger spaces we write the new index as the
last one and then for the $n$-simplex case we get the
$[\frac12n(n+1)+1]\times[\frac12n(n+1)+1]$ matrix 
\begin{equation}
\tilde A_{K_1\dots K_n}= \left[\begin{array}{cc}A_{K_1\dots
K_n}&B_{K_1\dots K_n} \cr 0&1\end{array}\right].
\end{equation}
[We use square brackets when writing out these index matrices.]  For
example the 2-simplex equation becomes
\begin{eqnarray}
&&\left[\begin{array}{cc}A_{12}&B_{12}\cr0&1\end{array}\right]
\left[\begin{array}{cc}A_{13}&B_{13}\cr0&1\end{array}\right]
\left[\begin{array}{cc}A_{23}&B_{23}\cr0&1\end{array}\right]\nonumber\\
&&\qquad=\left[\begin{array}{cc}A_{23}&B_{23}\cr0&1\end{array}\right]
\left[\begin{array}{cc}A_{13}&B_{13}\cr0&1\end{array}\right]
\left[\begin{array}{cc}A_{12}&B_{12}\cr0&1\end{array}\right],
\end{eqnarray}
and expanding this yields (\ref{E:Aeq}) for $A$ and
\begin{equation}
A_{12}A_{13}B_{23}+A_{12}B_{13}+B_{12}=
A_{23}A_{13}B_{12}+A_{23}B_{13}+B_{23},
\label{E:ABeq}
\end{equation}
for $B$. The higher simplex equations have equally simple matrix
form. In fact, formally the equations now look exactly as in
(\ref{E:2s}-\ref{E:4s}) with $\tilde A$ instead of $R$, but the
interpretation is different: for $\tilde A$ we have ordinary matrix
products.

\section{Symmetries}
Before starting to solve the equations it is necessary to discuss
their symmetries.  For one thing, we only want to list the basic
solutions from which the others are obtained by the allowed
transformations.  It is well known \cite{some} that the $n$-simplex
equations are form invariant under discrete transformations of index
transposition and index reversal. Now we should also see if these
transformations preserve the linear permutation structure and what
they imply on $A$ and $B$.

\subsection{Index transposition of $R$}
If $R_{i_1\dots i_n}^{j_1\dots j_n}$ is a solution of the $N$-simplex
equation, then $(IR)^{j_1\dots j_n}_{i_1\dots i_n}:=R^{i_1\dots
i_n}_{j_1\dots j_n}$ is also a solution.  This is easy to see from the
structure of the equation.

Let us now see what the above symmetry implies for the index matrix $A$.
From the definitions above it follows that
\begin{equation}
(IR)^{j_1\dots j_n}_{i_1\dots i_n}
=\prod_{\alpha=1}^n \delta^{i_\alpha}_{A_{\alpha}^{\beta}j_\beta+B_\alpha}
=\prod_{\alpha=1}^n \delta^{j_\alpha}_{(IA)_{\alpha}^{\beta}i_\beta+
(IB)_\alpha}.
\end{equation}
and by comparing the two expressions we find that if $R=\delta(A,B)$ is a
solution, then $(IR)=\delta(A^{-1},-A^{-1}B)$ is a solution, that is,
$(IA)=A^{-1}$, $(IB)=-A^{-1}B$.

It is easy to see that this is also an invariance of the $\tilde A$
equation.  As a matrix equation it is clearly invariant under matrix
inversion (which furthermore does not change the location of the
inserted pieces of the unit matrix) and $\left[\begin{smallmatrix}
A&B\\ 0&1\end{smallmatrix}\right]^{-1}=
\left[\begin{smallmatrix}A^{-1}&-A^{-1}B\\
0&1\end{smallmatrix}\right]$.

\subsection{Index reversal of $R$}
It is also easy to see that if $R_{i_1\dots i_n}^{j_1\dots j_n}$ is a
solution then $(CR)_{i_1\dots i_n}^{j_1\dots j_n}:=R_{i_n\dots
i_1}^{j_n\dots j_1}$ is a solution. We have
\begin{equation} 
(CR)^{j_1\dots j_n}_{i_1\dots i_n} =\prod_{\alpha=1}^n
\delta^{j_{\alpha}}_{
A_{n+1-\alpha}^{n+1-\beta}i_{\beta}+B_{n+1-\alpha}}
=\prod_{\alpha=1}^n \delta^{j_\alpha}_{(CA)_{\alpha}^{\beta}i_\beta+
(CB)_\alpha},  
\end{equation}
and the comparison yields $ (CA)_\alpha^\beta=
A_{n+1-\alpha}^{n+1-\beta},\, (CB)_\alpha=B_{n+1-\alpha}$, that is,
reflection across the center of the matrix or vector.

Since the $\tilde A$-equation is a matrix equation it is invariant
under any permutation of the set over which the summation is taken: if
$\tilde A_\alpha^\beta$ solves the equation, then $(\sigma \tilde
A)_\alpha^\beta:= \tilde A_{\sigma(\alpha)}^{\sigma(\beta)}$ where
$\sigma$ is any permutation operator, is also a solution.  However, we
have to keep intact the structure of inserted parts of the unit matrix
in the various terms, and then it appears that only the above reversal
is possible.

\subsection{Transposition of $A$}
From the point of view of $A$ the matrix equations have one more
discrete symmetry: they are invariant also under transposition.
However, this does not seem to correspond to any obvious invariance of
$R$. In the following it will turn out that often the transposition of
a solution $A$ is also obtained by the central reflection (accompanied
with parameter changes). However, this is not always true, and when it
is not, it turns out that often the accompanying $B$ will also be
different.

This is a rather interesting result from the point of view of studying
the structure of the equations. Normally imposing ansatze on the
solutions restrict the symmetries, because the symmetries of the
equation may not be symmetries of the ansatze. In the present case
this happens with the continuous transformation below. However, the
opposite can happen as well: in the present case the ansatz leads to a
new formulation which has its own obvious symmetries, and some of
these do not seem to have any counterpart at the original level.

\subsection{Gauge transformations}
The gauge transformation $R_{K_1\dots K_n} \to (QR)_{K_1\dots K_n}=
Q_{K_1}^{-1}\cdots Q_{K_n}^{-1}\,R_{K_1\dots K_n}\,Q_{K_1}\cdots
Q_{K_n}$, is also an invariance of the $n$-simplex equations.  Now
that $R$ is made out of delta-functions the transformation matrix $Q$
must also be of that form, i.e.,
\[
Q_i^j=\delta^j_{ui+v},\quad (Q^{-1})_j^k=\delta_{u^{-1}j-u^{-1}v}^k,\quad
u,v\in \Bbb Z_D.
\]
(If $D$ is not prime $u^{-1}$ is not always defined.)
A simple calculation shows that if $R=\delta(A,B)$ is a solution, then
$QR:=\delta(A,QB)$, where
\begin{equation}
(QB)_\alpha:=uB_\alpha+(1-\sum_{\gamma}A_{\alpha}^{\gamma})v,
\label{E:gauge}
\end{equation}
is also a solution.  Thus only $B$ can change, and we can in fact put
one $B_\beta=0$, if $\sum_{\gamma}A_{\beta}^{\gamma}\ne 1$. Later we
will find that for many solutions the inhomogeneous $B$ part is such
that it can be completely eliminated by this gauge transformation.

In order to understand this as an invariance of the equations
(\ref{E:Aeq},\ref{E:ABeq}) we note first that (\ref{E:Aeq}) can
be written as
\[
A_{12}A_{13}(1-A_{23})+A_{12}(1-A_{13})+(1-A_{12})=
A_{23}A_{13}(1-A_{12})+A_{23}(1-A_{13})+(1-A_{23}).  
\] 
If we now sum over the rightmost index of this equation and take its
linear combination with equation (\ref{E:ABeq}) we get (\ref{E:gauge})
for $(QB)$.

Matrix equations are invariant under a much larger group of similarity
transformations: $A\to O^{-1}AO$, but now that we have to preserve the
structure of having inserted pieces of the unit matrix in $A_{K_1\dots
K_n}$ these similarity transformation are allowed only with the matrix
$O=\left[\begin{smallmatrix}x\Bbb I & y\Bbb I\\ 0&1
\end{smallmatrix}\right]$ corresponding to the above.

\section{Results for the 2-simplex equation}
The details for the Yang-Baxter or 2-simplex case are as follows. In
the homogeneous case we write
\begin{equation}
R_{i_1i_2}^{j_1j_2}=\delta^{j_1}_{ai_1+bi_2}
\delta^{j_2}_{ci_1+di_2},
\end{equation}
so that $A= \left[\begin{smallmatrix} a&b\\
c&d\end{smallmatrix}\right]$ and then the 2-simplex equation becomes
\begin{equation}
\left[\begin{array}{ccc}a & b & 0\cr c & d & 0\cr 0 & 0 &
1\end{array}\right] \left[\begin{array}{ccc}a & 0 & b\cr 0 & 1 & 0\cr
c & 0 & d\end{array}\right] \left[\begin{array}{ccc}1 & 0 & 0\cr 0 & a
& b\cr 0 & c & d\end{array}\right] = \left[\begin{array}{ccc}1 & 0 &
0\cr 0 & a & b\cr 0 & c & d\end{array}\right]
\left[\begin{array}{ccc}a & 0 & b\cr 0 & 1 & 0\cr c & 0 &
d\end{array}\right] \left[\begin{array}{ccc}a & b & 0\cr c & d & 0\cr
0 & 0 & 1\end{array}\right].
\label{E:2sima}
\end{equation}
(In \cite{Korthe} similar matrices, but with the entries also being
matrices, are used to define a dynamical system.)  This yields 5
equations,
\begin{equation}
abc=0,\,
bcd=0,\,
bc(b-c)=0,\,
b(ad+b-1)=0,\,
c(ad+c-1)=0,
\end{equation}
whose solutions are discussed below.

\subsection{$D$ is prime}
Recall that we are working with integers modulo $D$, the number of
states.  If $D$ is prime there are no divisors of zero, and we can solve
the equations with conventional rules of algebra.  It is easy to show
that in this case there are precisely four nonsingular solutions:
\begin{equation}
A_2^{(1)}=\left[\begin{array}{cc}a&0\cr 0&d\end{array}\right],\,
A_2^{(2)}=\left[\begin{array}{cc}a&1-ad\cr 0&d\end{array}\right],\,
A_2^{(2r)}=\left[\begin{array}{cc}a&0\cr 1-ad&d\end{array}\right],\,
A_2^{(3)}=\left[\begin{array}{cc}0&1\cr 1&0\end{array}\right].
\end{equation}
Whether these are really different depends on $D$. Note for example,
that for $D=2$ cases (2) and (2r) reduce to case (1); the same holds with
$a=d=2$ when $D=3$. Solution (2r) is obtained from (2) by central
reflection, and there is actually no need to mention it separately.

When the inhomogeneous part $B=[x,y]^t$ is included we have to solve
(\ref{E:ABeq}), which amounts to
\begin{equation}
b(x+ay)=0,\quad c(y+dx)=0,\quad x(c+d-bc-1)=y(a+b-bc-1).
\end{equation}
The solutions then split further and we get
\begin{equation*}
[A|B]_2^{(1a)}=\left[\begin{array}{cc|c}1&0&x\cr
0&1&y\end{array}\right],\quad
[A|B]_2^{(1b)}=\left[\begin{array}{cc|c}a&0&(a-1)z\cr 0&d&(d-1)z
\end{array}\right],
\end{equation*}
\begin{equation*}
[A|B]_2^{(2)}=\left[\begin{array}{cc|c}a&1-ad&-az\cr
0&d&z\end{array}\right],\quad
[A|B]_2^{(3)}=\left[\begin{array}{cc|c}0&1&0\cr 1&0&0\end{array}\right].
\end{equation*}
However, we have not yet used the gauge freedom. For $[A|B]_2^{(1a)}$
and $[A|B]_2^{(3)}$ the row sums of $A$ are $=1$, and thus according to
(\ref{E:gauge}) we cannot change the inhomogeneous part, except by an
overall multiplication. For $[A|B]_2^{(1b)}$ the gauged inhomogeneous
part will turn out to be 
\[
\begin{bmatrix}(QB)_1\\(QB)_2\end{bmatrix}=
\begin{bmatrix}(a-1)(uz-v)\\(d-1)(uz-v)\end{bmatrix}
\]
and by choosing $v=uz$ we get $(QB)_i=0$. For
$[A|B]_2^{(2)}$ we get similarly 
\[
\begin{bmatrix}(QB)_1\\(QB)_2\end{bmatrix}=
\begin{bmatrix}-a(uz+v(1-d))\\uz+v(1-d)\end{bmatrix}
\]
Now if $d\ne 1$ we can again transform to $(QB)_i=0$, but if $d=1$
only scaling is possible. Thus the final form of the solutions of the
2-simplex case is
\begin{equation*}
[A|B]_2^{(1a)}=\left[\begin{array}{cc|c}1&0&x\cr 0&1&y\end{array}\right],\,
[A|B]_2^{(1b)}=\left[\begin{array}{cc|c}a&0&0\cr 0&d&0\end{array}\right],\,
[A|B]_2^{(3)}=\left[\begin{array}{cc|c}0&1&0\cr 1&0&0\end{array}\right],
\end{equation*}
\begin{equation*}
[A|B]_2^{(2a)}=\left[\begin{array}{cc|c}a&1-a&-az\cr
0&1&z\end{array}\right],\quad
[A|B]_2^{(2b)}=\left[\begin{array}{cc|c}a&1-ad&0\cr
0&d&0\end{array}\right],
\end{equation*}
up to the allowed transformations.

\subsection{$D=2,3$}
For $D=2$ the above yields basically two solutions for the 2-simplex
equation, $R$ is either the unit matrix (with possible
inhomogeneities) or the permutation matrix $P$,
\begin{equation}
[A|B]=
\left[\begin{array}{cc|c}
1 & 0 & x\\
0 & 1 & y\\
\end{array}\right]\mbox{ or }
\left[\begin{array}{cc|c}
0 & 1 & 0\\
1 & 0 & 0\\
\end{array}\right],
\label{E:d2}
\end{equation}
where $x,y\in\Bbb Z_2$. The same five solutions are obtained by a
brute force search without the linearity assumption.

The results (\ref{E:d2}) work for any $D$, but already for $D=3$
we get other homogeneous solutions, including triangular ones:
\begin{eqnarray*}
\left[\begin{array}{cc}1&0\cr 0&2\end{array}\right],\,
\left[\begin{array}{cc}2&0\cr 0&2\end{array}\right],\,
\left[\begin{array}{cc}1&2\cr 0&2\end{array}\right],\,
\left[\begin{array}{cc}2&2\cr 0&1\end{array}\right],\,
\left[\begin{array}{cc}2&2\cr 0&2\end{array}\right],
\end{eqnarray*} 
and their reflections.

\subsection{$D=4$}
The situation is quite different if $D=4$, because of divisors of
zero: $2\cdot 2=0 \pmod 4$. In addition to the above generic
solutions, we get new base solutions
\[
\left[\begin{array}{cc}1&2\cr 0&1\end{array}\right],\,
\left[\begin{array}{cc}1&2\cr 2&1\end{array}\right],\,
\left[\begin{array}{cc}1&2\cr 2&3\end{array}\right],\,
\left[\begin{array}{cc}3&2\cr 0&3\end{array}\right],\,
\left[\begin{array}{cc}3&2\cr 2&1\end{array}\right],\,
\left[\begin{array}{cc}3&2\cr 2&3\end{array}\right],\,
\]
and their reflections. When the inhomogeneous parts are added we get
\[
\left[\begin{array}{cc|c}1&2&0\cr 0&1&2x\end{array}\right],\,
\left[\begin{array}{cc|c}1&2&y+2x+gx\cr 2&1&y+gx\end{array}\right],\,
\left[\begin{array}{cc|c}1&2&0\cr 2&3&2x\end{array}\right],
\]\[
\left[\begin{array}{cc|c}3&2&2x\cr 0&3&0\end{array}\right],\,
\left[\begin{array}{cc|c}3&2&2x\cr 2&1&0\end{array}\right],\,
\left[\begin{array}{cc|c}3&2&y+2x\cr 2&3&y\end{array}\right].
\]
In the second solution there is no obvious way to fix the gauge
parameter $g$ and it has been left open.

Thus for $4$-state models there seem to be are additional symmetries
and solutions, and perhaps this case needs more detailed studies.

\section{Results for the 3-simplex equation}
Higher simplex equations have many reductions to lower simplex
equations, and it is not necessary to repeat them. For example, any
solutions of the 2-simplex equation generates a solution of the
3-simplex equation by $R_{ijk}=R_{ij}\delta_k$ or $\delta_i R_{jk}.$
These solutions (and those with $\det A=0$) will not be included in
the following list and the solution below are genuine 3-simplex
solutions.  Note also that $R_{ijk}=R_{ik}\delta_j$ (with $\delta$ on
the central index) is not automatically a solution, in particular the
permutation matrix $R_{i_1i_2i_3}^{j_1j_2j_3}= \delta_{i_1}^{j_3}
\delta_{i_2}^{j_2}\delta_{i_3}^{j_1}$ does not solve the tetrahedron
equation.

In order to solve the tetrahedron equation under the present ansatz we
first consider the homogeneous part. The equation to solve is just
like (\ref{E:3s}) with ${\cal R}$ replaced with $A$. When the matrix
\[
A=\left[\begin{array}{ccc}a & b & c\\x & y & z\\u & v & w
\end{array}\right]
\]
is inserted into the $6\times 6$ matrix $\tilde A_{K_1K_2K_3}$ the six
different ways indicated in (\ref{E:3s}) and we compute the
corresponding matrix product we find 29 equations:
\begin{eqnarray*}
&  a b x=0,\,
b x y=0,\,
v y z=0,\,
v w z=0,\,
&\\& 
b x (b - x)=0,\,
v z (v - z)=0,\,
y (b u - c v)=0,\,
y ( - c x + u z)=0,\,
&\\& 
b (a y + b - 1)=0,\,
x (a y + x - 1)=0,\,
z (w y + z - 1)=0,\,
v (w y + v - 1)=0,\,
&\\& 
a b u z + a c x + b c u=0,\,
b v x z + c v y + c x y=0,\,
b u w z + c u z + c v w=0,\,
&\\& 
 a b u + a c v x + c u x=0,\,
 b u y + b v x z + u y z=0,\,
 c u v + c v w x + u w z=0,\,
&\\& 
a b w z + a c z + b c w + c^2=0,\,
a u v + a v w x + u^2 + u w x=0,\,
&\\& 
b u v z + c u y + c v^2 - c v z=0,\,
b u x z - b c x + c u y + c x^2=0,\,
&\\& 
 - b^2 u - b c v x + b u x - c u y=0,\,
 - c u y - c v x z + u v z - u z^2=0,\,
&\\& 
b w x z + c w y + c x z + c z - c=0,\,
a b v z + a c y + b c v + b c - c=0,\,
&\\& 
 - a u y - a v x z - u x z - u x + u=0,\,
 - b u v - b v w x - u v - u w y + u=0,\,
&\\&
 - b c u + b u^2 z - c^2 v x + c u v + c u x - c u z=0.&
\end{eqnarray*}
By just considering the first four equations the problem can be split
into 9 different cases, and each one of them can then be solved rather
easily. After eliminating those solutions that reduce to 2-simplex
solutions and those with noninvertible $A$ we find 3 basic solutions
from which others are obtained by the allowed transformations. These
solutions and their nonhomogeneous additions will be discussed below.

\subsection{}
The first base solution is
\begin{equation*}A^{(1)}_3=
\left[\begin{array}{ccc}
0 & 1 & -d \\
1 & 0 & 1 \\
0 & 0 & d 
\end{array}\right]
\end{equation*}
and when inhomogeneities are added it splits into two:
\begin{equation*}[A|B]^{(1a)}_3=
\left[\begin{array}{ccc|c}
0 & 1 & -1 & x\\
1 & 0 & 1 & y\\
0 & 0 & 1 & 0
\end{array}\right],\quad [A|B]^{(1b)}_3=
\left[\begin{array}{ccc|c}
0 & 1 & -d & 0\\
1 & 0 & 1 & 0\\
0 & 0 & d &0
\end{array}\right]
\end{equation*}
For $[A|B]^{(1a)}_3$ $x$ or $y$ can be still eliminated by a gauge
transformation, for $[A|B]^{(1b)}_3$ the gauge freedom has already
been used above. 

The transpose of $A^{(1)}_3$ is not obtained by central reflection
and therefore constitutes another solution: 
\begin{equation*}[A|B]^{(1ta)}_3=
\left[\begin{array}{ccc|c}
0 & 1 & 0 & 0\\
1 & 0 & 0 & 0\\
-1 & 1 & 1 & z
\end{array}\right],\quad [A|B]^{(1tb)}_3=
\left[\begin{array}{ccc|c}
0 & 1 & 0 & 0\\
1 & 0 & 0 & 0\\
-d & 1 & d & 0
\end{array}\right]
\end{equation*}
These forms cannot be changed by gauge, except by $z\to uz$.

\subsection{}
There are two upper triangular solutions
\begin{equation*}A_3^{(2)}=\left[\begin{array}{ccc}
a & 1-ab & a(bc-1)\\
0 & b & 1-bc\\
0 & 0 & c
\end{array}\right],\quad A_3^{(2t)}=
\left[\begin{array}{ccc}
a & 1-ab & c(ba-1)\\
0 & b & 1-bc\\
0 & 0 & c
\end{array}\right].\end{equation*}
They differ only in the upper right hand entry and are related by
transposition and central reflection (followed by $a\leftrightarrow
c$). But since transposition is not a symmetry of the inhomogeneous
part they have to be analyzed separately.

Depending on which parameters have unit value we get three
solutions: 
\begin{equation*}[A|B]_3^{(2a)}=\left[\begin{array}{ccc|c}
1 & 1-b & b-1&x\\
0 & b & 1-b&-bz\\
0 & 0 & 1&z
\end{array}\right],\quad 
[A|B]_3^{(2b)}=\left[\begin{array}{ccc|c}
a & 1-b & a(b-1)&abz\\
0 & b & 1-b&-bz\\
0 & 0 & 1&z
\end{array}\right],\end{equation*}\begin{equation*}
[A|B]_3^{(2c)}=\left[\begin{array}{ccc|c}
a & 1-ab & a(bc-1)&0\\
0 & b & 1-bc&0\\
0 & 0 & c&0
\end{array}\right]. 
\end{equation*}
This solution illustrates nicely how added freedom in $A$ decreases
freedom in $B$.

For the transpose $A_3^{(2t)}$ we get two solutions
\begin{equation*}
[A|B]_3^{(2ta)}=\left[\begin{array}{ccc|c}
a & 1-a & a-1&x\\
0 & 1 & 0&y\\
0 & 0 & 1&0
\end{array}\right],\quad
[A|B]_3^{(2tb)}=\left[\begin{array}{ccc|c}
a & 1-ab & c(ba-1)&0\\
0 & b & 1-bc&0\\
0 & 0 & c&0
\end{array}\right]. 
\end{equation*}

\subsection{}
For the next solution the inhomogeneous terms can always be gauged
away and we have
\begin{equation*}[A|B]_3^{(3)}=
\left[\begin{array}{ccc| c}
a & 0 & 0&0\\
1-ab & b &1-bc&0\\
0 & 0 & c&0
\end{array}\right].\end{equation*}
Here we have assumed that at least one of $a,b,c$ is $\ne 1$, else we
get a diagonal solution with arbitrary $B$.

The transpose is again a separate case, we get first
\begin{equation*}[A|B]_3^{(3t)}=
\left[\begin{array}{ccc| c}
a & 1-ab & 0 & -ay\\
0 & b &    0 & y\\
0 & 1-bc & c & -cy
\end{array}\right].\end{equation*}
Now if $b\ne 1$ the inhomogeneous part can be eliminated, and we have
finally two solutions
\begin{equation*}[A|B]_3^{(3ta)}=
\left[\begin{array}{ccc| c}
a & 1-a & 0 & -ay\\
0 & 1 &    0 & y\\
0 & 1-c & c & -cy
\end{array}\right],\quad
[A|B]_3^{(3tb)}=
\left[\begin{array}{ccc| c}
a & 1-ab & 0 & 0\\
0 & b &    0 & 0\\
0 & 1-bc & c & 0
\end{array}\right]
\end{equation*}
Note how this case is built up from 2-simplex solutions $[A|B]^{(2b)}$,
but not as simple tensor products.

\subsection{$D=2$}
When all indices are modulo 2 only two solutions remain (in addition
to reducible ones) namely
\begin{equation*}
\left[\begin{array}{ccc|c}
0 & 1 & 1&0\\
1 & 0 & 1&y\\
0 & 0 & 1&0
\end{array}\right],\quad
\left[\begin{array}{ccc|c}
0 & 1 & 0 & 0\\
1 & 0 & 0 & 0\\
1 & 1 & 1 & y
\end{array}\right]
\end{equation*}
where $y\in\Bbb Z_2$ and we have used the gauge freedom to eliminate
$B_1$ in the first case. Here it might be useful to record the
corresponding $R$-matrices for $y=0$:
 
\[\left(\begin{array}{cccc|cccc}
1 & . & . & . & . & . & . & .\\
. & . & 1 & . & . & . & . & .\\
. & 1 & . & . & . & . & . & .\\
. & . & . & 1 & . & . & . & .\\
\hline
. & . & . & . & . & . & . & 1\\
. & . & . & . & . & 1 & . & .\\
. & . & . & . & . & . & 1 & .\\
. & . & . & . & 1 & . & . & .
\end{array}\right),
\quad\left(\begin{array}{cccc|cccc}
1 & . & . & . & . & . & . & .\\
. & . & . & . & . & . & 1 & .\\
. & . & . & . & . & 1 & . & .\\
. & . & . & 1 & . & . & . & .\\
\hline
. & . & . & . & 1 & . & . & .\\
. & . & 1 & . & . & . & . & .\\
. & 1 & . & . & . & . & . & .\\
. & . & . & . & . & . & . & 1
\end{array}\right).
\]
These bear some resemblance with know solutions \cite{Kor,Label,Ser,some}.

\section{Results for the 4-simplex equation}
For the 4-simplex case the $4\times4$ index matrix is embedded into
$10\times10$ matrices in four ways.  The equations resulting from
(\ref{E:4s}) were solved using the groebner-package of REDUCE
\cite{Gro}.  From the results we eliminated those solutions for which
$A$ was in a block form corresponding to $R$'s with tensor products
form $R_{i_1i_2i_3i_4}^{j_1j_2j_3j_4}=
\delta_{i_1}^{aj_1}M_{i_2i_3i_4}^{j_2j_3j_4}$,
$R_{i_1i_2i_3i_4}^{j_1j_2j_3j_4}=
M_{i_1i_2i_3}^{j_1j_2j_3}\delta_{i_4}^{aj_4}$ or
$R_{i_1i_2i_3i_4}^{j_1j_2j_3j_4}=
K_{i_1i_2}^{j_1j_2}L_{i_3i_4}^{j_3j_4}$ where $M$ is a solution of the
3-simplex equation and $K,L$ of the 2-simplex equation.  From the
remaining list we eliminated all cases obtained from the basic ones by
central reflection or by inverse, and those with singular $A$.
Furthermore we considered only the generic case of a prime $D$.

The solutions $A$ of the homogeneous equation (and their transposes)
were next used as starting points for constructing the non-homogeneous
part $B$.  Then the continuous gauge freedom was applied to eliminate
some freedom from $B$.  The final result is as follows:

\subsection*{Permutation blocks}
\subsection{}
\begin{equation*}\left[\begin{array}{cccc|c}
0 & 1 & 0 & -1&b_1\\
1 & 0 & 1 & 0&b_2\\
0 & 0 & 0 & 1&0\\
0 & 0 & 1 & 0&0
\end{array}\right]\end{equation*}
After a gauge transformation (\ref{E:gauge}) we would get
$B^t=[b_1+v,b_2-v,0,0]$ and we could eliminate either $b_1$ or $b_2$.

\subsection{}
\begin{equation*}\left[\begin{array}{cccc|c}
0 & 1 & 0 & 0&0\\
1 & 0 & 0 & 0&0\\
 - a & 1 & a & 1-ab&0\\
 0 & 0 & 0 & b&0 
\end{array}\right]\end{equation*}
In the generic case we get $B^t=[0,0,z(b-1),-z(ab-1)]$ but this
can be eliminated by the gauge transformation.  Only if $a=b=1$ would we
get something that cannot be gauged away, but in that case the system
reduces to a 3-simplex solution. 

\subsection{}
The transpose of the above solutions is a separate case, and yields
\begin{equation*}
\left[\begin{array}{cccc|c}
0 & 1 & -1 & 0&b_1\\
1 & 0 & 1 & 0&b_2\\
 0 & 0 & 1 & 0&0\\
 0 & 0 & 1-b & b&0 
\end{array}\right],\,
\left[\begin{array}{cccc|c}
0 & 1 & -a & 0&0\\
1 & 0 & 1 & 0&0\\
 0 & 0 & a & 0&0\\
 0 & 0 & 1-ab & b&0 
\end{array}\right],
\end{equation*}
(where $b_1$ or $b_2$ could be gauged away). Again when the $A$ part
is restricted the $B$ part gains some freedom.

\subsection{}
The next cases are somewhat similar to the above, we get
\begin{equation*}
\left[\begin{array}{cccc|c}
0 & 1 & -a & a-1&b_1\\
1 & 0 & 1 & 0&b_2\\
 0 & 0 & a & 1-a&0\\
 0 & 0 & 0 & 1&0 
\end{array}\right],\,
\left[\begin{array}{cccc|c}
0 & 1 & -a & ab-1&0\\
1 & 0 & 1 & 0&0\\
 0 & 0 & a & 1-ab&0\\
 0 & 0 & 0& b&0 
\end{array}\right].
\end{equation*}

\subsection{}
\begin{equation*}
\left[\begin{array}{cccc|c}
0 & 1 & 0 & 0&0\\
1 & 0 & 0 & 0&0\\
-a & 1 & a & 0&0\\
a-1 & 0 & 1-a& 1&x 
\end{array}\right],\,
\left[\begin{array}{cccc|c}
0 & 1 & 0 & 0&0\\
1 & 0 & 0 & 0&0\\
-a & 1 & a & 0&0\\
ab-1 & 0 & 1-ab& b&0 
\end{array}\right]
\end{equation*}

\subsection{}
\begin{equation*}
\left[\begin{array}{cccc|c}
0 & 1 &  - a & a&ax\\
1 & 0 & 1 &  - 1&-x\\
0 & 0 & a & 1-a&-ax\\
0 & 0 & 0 & 1&x
\end{array}\right],\,
\left[\begin{array}{cccc|c}
0 & 1 &  - a & a b&0\\
1 & 0 & 1 &  - b&0\\
0 & 0 & a & 1-ab&0\\
0 & 0 & 0 & b&0
\end{array}\right].
\end{equation*}

\subsection{}
For the transpose of the above nothing can be gauged away, 
and we get
\begin{equation*}
\left[\begin{array}{cccc|c}
0 & 1 &  0 & 0&0\\
1 & 0 & 0 &  0&0\\
-1 & 1 & 1 & 0&x\\
1 & -1 & 0 & 1&y
\end{array}\right],\,
\left[\begin{array}{cccc|c}
0 & 1 &  0 & 0&0\\
1 & 0 & 0 &  0&0\\
-1 & 1 & 1 & 0&-bx\\
b & -b & 1-b & b&x
\end{array}\right],
\end{equation*}
\begin{equation*}
\left[\begin{array}{cccc|c}
0 & 1 &  0 & 0&0\\
1 & 0 & 0 &  0&0\\
-a & 1 & a & 0&0\\
a & -1 & 1-a & 1&x
\end{array}\right],\,
\left[\begin{array}{cccc|c}
0 & 1 &  0 & 0&0\\
1 & 0 & 0 &  0&0\\
-a & 1 & a & 0&0\\
ab & -b & 1-ab & b&0
\end{array}\right].
\end{equation*}

\subsection{}
The next $A$ matrix is invariant under central reflection, and gauge
transformation changes nothing. We get three different cases
\begin{equation*}
\left[\begin{array}{cccc|c}
1 & 1 & -1 & 0 & x\\
0 & 0 & 1 & 0 & 0\\
0 & 1 & 0 & 0 & 0\\
0 & -1 & 1 & 1 & y
\end{array}\right],\,
\left[\begin{array}{cccc|c}
1 & 1 & -1 & 0 & 0\\
0 & 0 & 1 & 0 & 0\\
0 & 1 & 0 & 0 & 0\\
0 & -d & 1 & d & y
\end{array}\right],\,
\left[\begin{array}{cccc|c}
a & 1 & -a & 0 & 0\\
0 & 0 & 1 & 0 & 0\\
0 & 1 & 0 & 0 & 0\\
0 & -d & 1 & d & 0
\end{array}\right].
\end{equation*}

\subsection{}
For the transpose of the above the inhomogeneous part is quite
different and we get
\begin{equation*}
\left[\begin{array}{cccc|c}
1 & 0 & 0 & 0 & x\\
1 & 0 & 1 & -1 & y\\
-1 & 1 & 0 & 1 & -x\\
0 & 0 & 0 & 1 & -y
\end{array}\right],\,
\left[\begin{array}{cccc|c}
a & 0 & 0 & 0 & 0\\
1 & 0 & 1 & -d & 0\\
-a & 1 & 0 & 1 & 0\\
0 & 0 & 0 & d & 0
\end{array}\right].
\end{equation*}
In both cases there are two free parameters.

\subsection{}
\begin{equation*}
\left[\begin{array}{cccc|c}
a & 1 & -a & a & 0\\
0 & 0 & 1 & -1 & 0\\
0 & 1 & 0 & 1 & y\\
0 & 0 & 0 & 1 & 0
\end{array}\right],\,
\left[\begin{array}{cccc|c}
a & 1 & -a & ad & 0\\
0 & 0 & 1 & -d & 0\\
0 & 1 & 0 & 1 & 0\\
0 & 0 & 0 & d & 0
\end{array}\right]
\end{equation*}

\vskip 1cm\noindent{\large{\bf {Triangular blocks}}}
\subsection{}
\begin{equation*}
\left[\begin{array}{cccc|c}
1 & 1-b & b - 1 & 1-b & 0 \\
0 & b & 1-b & b - 1 & 0 \\
0 & 0 & 1 &  0 & x \\
0 & 0 & 0 & 1 & 0
\end{array}\right],\,
\left[\begin{array}{cccc|c}
a & 1-ab & c (a b - 1) & c d (1-ab) & 0 \\
0 & b & 1-bc & d (b c - 1) & 0 \\
0 & 0 & c &  1-cd & 0 \\
0 & 0 & 0 & d & 0
\end{array}\right].
\end{equation*}

\subsection{}
\label{6.12}
\begin{equation*}
\left[\begin{array}{cccc|c}
1 & 0 & 0 & 0 & x \\
1-b & b & 0 & 0 & -bx \\
b-1 & 1-b & 1 & 0 & y \\
1-b & b-1 & 0 & 1 & z
\end{array}\right],\,
\left[\begin{array}{cccc|c}
1 & 0 & 0 & 0 & x \\
1-b & b & 0 & 0 & -bx \\
b-1 & 1-b & 1 & 0 & y \\
d(1-b) & d(b-1) & 1-d & d & -dy
\end{array}\right],
\end{equation*}
\begin{equation*}
\left[\begin{array}{cccc|c}
1 & 0 & 0 & 0 & x \\
1-b & b & 0 & 0 & -xb \\
c(b-1) & 1-bc & c & 0 & cbx \\
c(1-b) & (bc-1) & 1-c & 1 & y
\end{array}\right],\,
\left[\begin{array}{cccc|c}
a & 0 & 0 & 0 & 0 \\
1-ab & b & 0 & 0 & 0 \\
c(ab-1) & 1-bc & c & 0 & 0 \\
cd(1-ab) & d(bc-1) & 1-cd & d & 0
\end{array}\right].
\end{equation*}

\subsection{}
The detailed analysis of this case leads to some subcases that are
identical to those of central reflected (\ref{6.12}) and are not
repeated here.
\begin{equation*}
\left[\begin{array}{cccc|c}
a & 1-ab & a (b - 1) & a (1-b) & 0 \\
0 & b & 1-b & (b - 1) & 0 \\
0 & 0 & 1 & 0 & y \\
0 & 0 & 0 & 1 & 0
\end{array}\right],\,
\left[\begin{array}{cccc|c}
a & 1-ab & a (b c - 1) & a d (1-bc) & 0 \\
0 & b & 1-bc & d (b c - 1) & 0 \\
0 & 0 & c &  1-cd & 0 \\
0 & 0 & 0 & d & 0
\end{array}\right].
\end{equation*}
 
\subsection{}
\label{6.14}
\begin{equation*}
\left[\begin{array}{cccc|c}
a & 1-ab & 0 & 0 & 0 \\
0 & b & 0 & 0 & 0 \\
0 & 1-bc & c &  1-cd & 0 \\
0 & 0 & 0 & d & 0
\end{array}\right]
\end{equation*}

\subsection{}
\begin{equation*}
\left[\begin{array}{cccc|c}
1 & 1-b & b - 1 & 0 & x \\
0 & b & 1-b & 0 & -bz \\
0 & 0 & 1 &  0 & z \\
0 & 0 & 1-d & d & -dz
\end{array}\right],\,
\left[\begin{array}{cccc|c}
a & 1-a & a - 1 & 0 & 0 \\
0 & 1 & 0 & 0 & y \\
0 & 0 & 1 &  0 & 0 \\
0 & 0 & 1-d & d & 0
\end{array}\right],
\end{equation*}
\begin{equation*}
\left[\begin{array}{cccc|c}
a & 1-ab & c (a b - 1) & 0 & 0 \\
0 & b & 1-bc & 0 & 0 \\
0 & 0 & c &  0 & 0 \\
0 & 0 & 1-cd & d & 0
\end{array}\right].
\end{equation*}

\subsection{}
\begin{equation*}
\left[\begin{array}{cccc|c}
a 	& 0 	& 0 & 0 & 0 \\
1-ab 	& b 	& 0 & 0 & 0 \\
c(ab-1) & 1-bc 	& c & 1-cd & 0 \\
0 	& 0 	& 0 & d & 0
\end{array}\right]
\end{equation*}

\subsection{}
\begin{equation*}
\left[\begin{array}{cccc|c}
a & 1-ab & a (b- 1) &0 & abx \\
0 & b & 1-b & 0 & -bx \\
0 & 0 & 1 & 0 & x \\
0 & 0 &  1-d & d & -dx
\end{array}\right],\,
\left[\begin{array}{cccc|c}
a & 1-ab & a (b c - 1) &0 & 0 \\
0 & b & 1-bc & 0 & 0 \\
0 & 0 & c & 0 & 0 \\
0 & 0 &  1-cd & d & 0
\end{array}\right]
\end{equation*}

\subsection{}
\begin{equation*}
\left[\begin{array}{cccc|c}
1 	& 0 	& 0 & 0	 & x \\
1-b 	& b 	& 0 & 0	 & -xb \\
bc-1	& 1-bc 	& c &1-c & -x(1-bc)-cy \\
0 	& 0 	& 0 & 1	 & y
\end{array}\right],\,
\left[\begin{array}{cccc|c}
a 	& 0 	& 0 & 0 & 0 \\
1-ab 	& b 	& 0 & 0 & 0 \\
a(bc-1) & 1-bc 	& c & 1-cd & 0 \\
0 	& 0 	& 0 & d & 0
\end{array}\right]
\end{equation*}
 
\subsection{}
Here and in the following case we have rational entries in the index
matrix.
\begin{equation*}
\left[\begin{array}{cccc|c}
a & 1-a & 0 & 0 & 0 \\
0 & 1 & 0 & 0 & 0 \\
0 & 0 & 1 & 0 & x \\
0 & d-1 & 1-d & d & 0
\end{array}\right],\,
\left[\begin{array}{cccc|c}
a & 1-ab & 0 & 0 & 0 \\
0 & b & 0 & 0 & 0 \\
0 & 0 & 1/b & 0 & 0 \\
0 & d-b & 1-d/b & d & 0
\end{array}\right]
\end{equation*}

\subsection{}
\begin{equation*}
\left[\begin{array}{cccc|c}
a & 0 & 0 & 0 & 0 \\
1-ab & b & 0 & d-b & 0 \\
0 & 0 & 1/b & (b-d)/b & 0 \\
0 & 0 & 0 & d & 0
\end{array}\right]
\end{equation*}

\subsection{$D=2$}
For $D=2$ we have the following new solutions
\begin{equation*}
\left[\begin{array}{cccc|c}
0 & 1 & 0 & 1 & 0 \\
1 & 0 & 1 & 0 & 0 \\
0 & 0 & 0 & 1 & 0 \\
0 & 0 & 1 & 0 & 0
\end{array}\right],\,
\left[\begin{array}{cccc|c}
0 & 1 & 1 & 1 & 0 \\
1 & 0 & 1 & 1 & 0 \\
0 & 0 & 1 & 0 & 0 \\
0 & 0 & 0 & 1 & 0
\end{array}\right],
\end{equation*}
\begin{equation*}
\left[\begin{array}{cccc|c}
1 & 1 & 1 &0 & 0 \\
0 & 0 & 1 & 0 & 0 \\
0 & 1 & 0 & 0 & 0 \\
0 & 1 & 1 & 1 & 0
\end{array}\right],\,
\left[\begin{array}{cccc|c}
1 & 1 & 1 & 1 & 0 \\
0 & 0 & 1 & 1 & 0 \\
0 & 1 & 0 & 1 & 0 \\
0 & 0 & 0 & 1 & 0
\end{array}\right].
\end{equation*}

\section{Discussion}
In presenting this (complete) set of linear permutation type solutions
we hope that some of them could be used in other studies. These
applications may require further conditions, but we believe that
permutation type solutions are so benign that they should satisfy
these conditions, if just independence of any spectral parameter is
acceptable.

Another hope is that the solutions can teach us something about the
equations themselves.  Oner observation in that direction is that some
of the solutions fall into patterns that seem to continue to any $n$.
For example, $A_2^{(2)}$, $A_3^{(3t)}$ and $\ref{6.14}$ start a
pattern that seems to continue as
\[
\begin{bmatrix}
a_1 & 1-a_1a_2 & 0 & 0 & 0 &\dots \\
0   &	a_2    & 0 & 0 & 0 &\dots \\
0   & 1-a_2a_3 & a_3 & 1- a_3a_4 & 0 & \dots\\
0 & 0 & 0 & a_4 & 0 &  \\
0 & 0 & 0  & 1-a_4a_5 & a_5 &\ddots\\
\vdots & \vdots & \vdots &  &\ddots & \ddots
\end{bmatrix}.
\]
This band structure could make sense even as an infinite matrix,
and perhaps we should soon start to think what kind of object the
``$\infty$-simplex'' equation might be.

\subsection*{Acknowledgments}
I would like to thank F. Nijhoff and S. Sergeev for discussions and
for comments on the manuscript.

\end{document}